\documentclass[12pt]{amsart}

\def\nn{{\mathbb N}}
\def\tnpm{{\mathsf{t{\rm -}NM_k(s_1,\ldots, s_t)}}}
\def\nmkn{{\mathsf {NM}}_k(n)}
\def\bnmkrs{{\mathsf {BNM}}_k(r,s)}
\def\npmn{{\mathsf {NPM}}_{2k}}
\def\npm{{\mathsf {NPM}}}
\def\tt{{\mathsf {TT}}}
\def\ttk{{\mathsf {TT}}(|V|)}
\def\fc{{\mathsf {FC}}}
\def\bfc{{\mathsf {BFC}}}
\def\nfc{{\mathsf {NFC}}}
\def\con{{\mathrm {con}}}
\def\nbfc{{\mathsf {NBFC}}}
\def\P{{\mathsf P}}
\def\rh{\widetilde{H}}
\def\gp{G^\prime}
\def\gm{G^{-}}
\def\hp{H^\prime}

\def\A{{\mathcal A}}
\def\B{{\mathcal B}}
\def\C{{\mathcal C}}
\def\D{{\mathcal D}}
\def\M{{\mathcal M}}

\newtheorem{theorem}{Theorem}[section]
\newtheorem{lemma}[theorem]{Lemma}
\newtheorem{proposition}[theorem]{Proposition}
\newtheorem{corollary}[theorem]{Corollary}

\newtheorem{definition}[theorem]{Definition}
\newtheorem{remark}[theorem]{Remark}
\newtheorem{conjecture}[theorem]{Conjecture}
\newtheorem{question}[theorem]{Question}

\begin{document}

\title[Complexes of graphs]{Complexes of graphs with bounded matching size}
\date{\today}

\author[Linusson]{Svante Linusson$^1$}
\footnotetext[1]{Supported by EC's IHRP program through grant
HPRN-CT-2001-00272}
\address{Department of Mathematics\\ Lin\-k\"o\-ping University \\
SE-581 83 Link\"oping\\ Sweden} \email{linusson@mai.liu.se}

\author[Shareshian]{John Shareshian$^2$}
\footnotetext[2]{Partially supported by National Science
Foundation grants DMS-0070757 and DMS-0030483}
\address{Department of Mathematics\\
         Washington University\\
         St Louis, MO 63130\\ USA }
\email{shareshi@math.wustl.edu}
\author[Welker]{Volkmar Welker$^3$}
\footnotetext[3]{Supported by EC's IHRP program through grant
HPRN-CT-2001-00272}
\address{Philipps-Universit\"at Marburg \\
Fachbereich Mathematik und Informatik \\ 
D-35032 Marburg \\ Germany  }
\email{welker@mathematik.uni-marburg.de}


\begin{abstract}
For positive integers $k,n$, we investigate the simplicial complex
$\nmkn$ of all graphs $G$ on vertex set $[n]$
such that every matching in $G$ has size less than $k$.  This 
complex (along with other associated cell complexes) is found to be
homotopy equivalent to a wedge of spheres. The number and dimension of the
spheres in the wedge are determined, and (partially conjectural) links to other
combinatorially defined complexes are described.
In addition we study for positive integers $r,s$ and $k$ the simplicial
complex $\bnmkrs$ of all bipartite graphs $G$ on bipartition $[r] \cup [\bar{s}]$ 
such that there is no matching of size $k$ in $G$, and obtain results similar to
those obtained for $\nmkn$.
\end{abstract}

\maketitle

\section{Introduction} \label{intro}

A monotone graph property is a collection $\Delta$ of (simple, loopless)
graphs on a fixed labelled vertex set $V$ such that
\begin{itemize}
\item[(C)] 
if $G \in \Delta$ and $H$ can be obtained from $G$ by removing an
edge then $H \in \Delta$, and
\item[(P)] if $G \in \Delta$ and $H$ can be obtained from $G$ by relabeling
the
vertices then $H \in \Delta$.
\end{itemize}

Condition (C) allows one to associate to each monotone graph
property $\Delta$ an abstract simplicial complex (also called
$\Delta$) whose $k$-di\-mensio\-nal faces are (indexed by) those
graphs in $\Delta$ which have $k+1$ edges, and condition (P) then
says that the natural action of the symmetric group $S_V$ on $V$
determines a simplicial action of $S_V$ on $\Delta$.  These facts
were used in the paper \cite{KSS} to apply results from algebraic
topology in attacking the evasiveness conjecture.  More recently,
certain classes of monotone graph properties arose naturally in
various problems from topology, algebra and combinatorics (see
\cite{BBLSW,Wa} for references).  

If $V = V_1 \cup V_2$ is a non-trivial partition of $V$ then one can consider
the following variant of condition (P):

\begin{itemize}
\item[(P')] each $G \in \Delta$ is bipartite with fixed bipartition $V_1
\cup V_2$.
If $G \in \Delta$ and $H$ can be obtained from $G$ by relabeling the
vertices within $V_1$ and within $V_2$ then $H \in \Delta$.
\end{itemize}

Graph properties satisfying (C) and (P') admit a natural action of
$S_{V_1} \times S_{V_2}$ and have also appeared recently
\cite{BBLSW,Wa}. Moreover, in this situation the evasiveness conjecture 
actually has been solved \cite{Y}.

Motivated only by curiosity, we
have studied the following two complexes:

\begin{itemize}
\item[$\rightarrow$] Let $n$ and $k$ be positive integers.
The complex $\nmkn$ consists of all graphs on vertex set $[n] := \{ 1,
\ldots, n\}$ 
which contain no matching of size $k$.
\item[$\rightarrow$] Let $r$, $s$ and $k$ be positive integers.
The complex $\bnmkrs$ consist of all bipartite graphs with
bipartite vertex set $[r] \cup [\bar{s}] := \{ 1, \ldots, r,
\bar{1},\ldots, \bar{s}\}$ 
which do not contain a matching of size $k$.
\end{itemize}

These complexes turn out to have very nice
topological structure which is (mysteriously and partially
conjecturally) related to that of some other combinatorially
defined complexes, as will be described below. Our proofs make use
of the discrete Morse theory of Forman (see \cite{Fo}).  This
theory has been applied several times in the study of monotone
graph properties (see \cite{BBLSW,Jo,LiSh,Sh}). However, our
results seem to be the first which make use of nonrudimentary
results in graph theory. Namely, we use the Gallai-Edmonds
structure theorem (see for example \cite{LoPl}).

Our main results are as follows.

\begin{theorem}
For $n,k \in \nn$, let $\Pi_{n-1}^1(k)$ be the set of all partitions
$\tau$ of $[n-1]$ into $n-2k+1$ subsets $\tau_1,\ldots,\tau_{n-2k+1}$ of
odd size.  Then $\nmkn$ has
the homotopy type of a wedge of spheres of dimension $3k-4$.  The number of
spheres in
this wedge is
\[
\sum_{\tau \in \Pi_{n-1}^1(k)}(\prod_{i=1}^{n-2k+1}
(\tau_i-2)!!)^2.
\]
\label{main-general}
\end{theorem}

A special case of this theorem is the following.

\begin{corollary}
For $k \in \nn$, let $\npmn$ be the complex of all graphs on vertex set
$[2k]$ which have
no perfect matching.  Then
\[
\npmn \simeq \bigvee_{(2k-3)!!^2}S^{3k-4}.
\]
\label{noper}
\end{corollary}

Calculating the (reduced) Euler characteristic of $\nmkn$ in two ways gives
the following
enumerative result.

\begin{corollary}
For $n,k \in \nn$, we have
\[
\sum_{G \in \nmkn}(-1)^{|E(G)|}=(-1)^{k-1}\sum_{\tau \in \Pi_{n-1}^k}
(\prod_{i=1}^{n-2k+1}(\tau_i-2)!!)^2.
\]
In particular,
\[
\sum_{G \in \npmn}(-1)^{|E(G)|}=(-1)^{k-1}(2k-3)!!^2.
\]
\label{enum}
\end{corollary}

\begin{theorem} \label{main-bipartite}
For $r,s,k \in \nn$, 
the homotopy type of $\bnmkrs$ is a wedge of spheres of dimension $2k-3$.
The number of spheres in
this wedge is $\binom{r-1}{k-1} \binom{s-1}{k-1}$.
\end{theorem}

Calculating the (reduced) Euler characteristic of $\bnmkrs$ in two ways
gives the following
enumerative result.

\begin{corollary}
For $r,s,k \in \nn$, we have
\[
\sum_{G \in \bnmkrs}(-1)^{|E(G)|}= \binom{r-1}{k-1} \binom{s-1}{k-1}.
\]
\label{enumbipartite}
\end{corollary}

In order to prove Theorem \ref{main-general}, 
we are forced to examine a certain quotient CW-complex.
Let $n=2m-1$ be odd.  A
graph $G$ on vertex set $[n]$ is called factor critical if
for each vertex $v$ of $G$, the graph obtained from $G$ by
removing $v$ and all edges containing $v$ has a perfect matching.
The classification of bipartite factor critical graphs is very simple.

\begin{remark} \label{critical-bipartite} 
A factor critical bipartite graph 
consists of a single vertex. 
\end{remark}

The set $\nfc_n$ of not factor critical graphs on $[n]$ is a
monotone graph property and therefore a subcomplex of the simplex
$\Sigma(n)$ whose vertices are the $\binom{n}{2}$ edges 
$\binom{[n]}{2} := \{ \{ i,j\} ~|~1 \leq i < j \leq n\}$ of the 
complete graph on vertex set $[n]$.  
The complex $\fc_n$ of factor critical
graphs is the quotient space $\Sigma(n)/\nfc_n$.  This space admits
an obvious cell decomposition, as described in Section \ref{dmt}.

\begin{theorem}
Let $n=2m-1 \in \nn$.  Then $\fc_n$ has the homotopy type of a
wedge of $(n-2)!!^2$ spheres of dimension $3m-4$. \label{fcth}
\end{theorem}

If $\Gamma$ is a nonempty subcomplex of the simplex $\Sigma$ then
$\Sigma/\Gamma$ is homotopy equivalent to the suspension of
$\Gamma$.  (One can prove this by showing that the mapping cone
of the identity embedding of $\Gamma$ into $\Sigma$ is the union of
two contractible subspaces whose intersection is $\Gamma$.) This allows
us to deduce the following result.

\begin{corollary}
Let $n=2m-1 \in \nn$.  Then $\nfc_n$ has the homotopy type of a
wedge of $(n-2)!!^2$ spheres of dimension $3m-5$.  \label{nfcth}
\end{corollary}

\begin{proof}
Our claim can be confirmed by direct observation when $n<5$, so
assume that $n \geq 5$.  We have $\fc_n=\Sigma(n)/\nfc_n$, and it
follows from Theorem \ref{fcth} and the following remark that
$\nfc_n$ has the same homology as the given wedge of spheres.
Since no graph on $[n]$ with three edges is factor critical, the
complex $\nfc_n$ contains the entire $2$-skeleton of $\Sigma(n)$
and is therefore simply connected. The corollary now follows from
the uniqueness of Moore spaces (see for example \cite[p.
368]{Hat}).
\end{proof}

By Remark \ref{critical-bipartite} the only factor critical bipartite 
graphs are singletons. Thus for our purposes it is useful to replace 
the concept of factor critical graphs by a new concept in the bipartite 
setting. We say a bipartite graph with bipartition $X \cup Y$ is
$q$-factor critical if $|X| = q$, $|Y| > q$ and for each $y \in Y$ the
graph $G-y$ has a matching of size $q$. The set $\nbfc(q,s)$ of bipartite
graphs on 
vertex set $[q] \cup [\overline{s}]$ which are not $q$-factor critical form a
subcomplex of the simplex $\Sigma(q,s)$ whose vertices are the
edges in the complete bipartite graph on bipartition $[q] \cup
[\overline{s}]$. (Here
$[\overline{s}]:=\{\overline{1},\ldots,\overline{s}\}$.)
The complex $\bfc(q,s)$ of $q$-factor critical 
bipartite graphs on $[q] \cup [\overline{s}]$ is the quotient complex 
$\Sigma(q,s) / \nbfc(q,s)$. 
Once again, this complex admits an obvious
cell decomposition (see Section \ref{dmt}).   
We will prove the following result, from which Theorem \ref{main-bipartite} 
will follow after some straightforward additional arguments.

\begin{theorem} \label{bfcth}
For $1 \leq q < s$ the complex $\bfc(q,s)$ is homotopy equivalent to a wedge
of $\binom{s-1}{q}$ spheres of dimension $2q-1$.
\end{theorem}

As in the non-bipartite case we deduce the following corollary.

\begin{corollary} For $1 \leq q <s $ the complex $\nbfc(q,s)$ is homotopy
equivalent to a wedge of $\binom{s-1}{q}$ spheres of dimension $2q-2$.
\end{corollary}

All of the complexes $\npm_n$, $\fc_n$ and $\nfc_n$ have homology
concentrated in a single dimension, and in each case the rank of
the unique nontrivial homology group is $(n-3)!!^2$ if $n$ is even
and $(n-2)!!^2$ if $n$ is odd. Moreover, in each case the action
of $S_n$ on the vertex set $[n]$ determines a simplicial (or
cellular) action of $S_n$ on the given complex. This action
determines a representation of $S_n$ on the nontrivial homology
group.  Representations of the groups $S_n$ whose degrees are the
ranks of these homology groups have arisen previously in work of
Calderbank, Hanlon and Robinson (\cite{CaHaRo}), Hanlon and Wachs
(\cite{HaWa}), and Hanlon (\cite{Ha}).  Fix $k \in \nn$ and for
each $n \in \nn$ let $\Pi_n^{1,k}$ be the subposet of the
partition lattice $\Pi_n$ consisting of those nontrivial proper
partitions in which each part has size equal to $1 \bmod k$.  In
\cite{Bj}, Bj\"orner showed that the order complex
$\Delta\Pi_n^{1,k}$ is Cohen-Macaulay and therefore has a unique
nontrivial homology group, which occurs in the top dimension
$t:=\lfloor \frac{n-2}{k} \rfloor-1$.  In \cite{HaWa}, it is shown
that character realized by the representation of $S_n$ on $\rh_t(\Delta\Pi_n^{1,k})$
is equal to the character 
${\sf lie}^{(k)}_n$ realized by the action of $S_n$ on a particular 
subspace of a vector space with a $(k+1)$-ary operation, called a 
Lie-$k$-algebra.  In \cite{CaHaRo}, the
character ${\sf lie}^{(2)}_n$ is determined. It is shown that if $n$ is odd then this character
has degree $(n-2)!!^2$.  
In \cite{Ha}, when $n \equiv 2 \bmod k$ a poset ${\mathcal
L}^{(k)}_n$ whose elements are trees with $n$ leaves labelled
bijectively with $[n]$ and nonleaves having degree equal to $2
\bmod k$ (but not equal to $2$) is defined.  It is shown that this
poset is Cohen-Macaulay of dimension $t$ and that the character
for the representation of $S_n$ on $\rh_t(\Delta{\mathcal L}^{(k)}_n)$ is
\[
{\sf lie}^{(k)}_{n-1} \uparrow_{S_{n-1}}^{S_n} - {\sf
lie}^{(k)}_n.
\]
Here $\uparrow$ denotes induction of representations or characters, and below
$\downarrow$ will represent restriction. Moreover, the dimension
of $\rh_t(\Delta{\mathcal L}^{(2)}_n)$ is $(n-3)!!^2$.  Let
$\varepsilon_n$ be the sign character of $S_n$. Wachs and the second
author of this paper have proved the following result

\begin{theorem}[Shareshian-Wachs]
Let $n=2k \in \nn$ be even.  Then
\[
\rh_{3k-4}(\npm_n)\downarrow^{S_n}_{S_{n-1}} \cong_{S_{n-1}}
\rh_{3k-5}(\nfc_{n-1}) \cong_{S_{n-1}}
\rh_t(\Delta\Pi_{n-1}^{1,2}) \otimes \varepsilon_n.
\]
\label{conj1}
\end{theorem}

We conjecture the following stronger result.

\begin{conjecture}
Let $n=2k$ be even.  Then
\[
\rh_{3k-4}(\npm_n) \cong_{S_n} \rh_t(\Delta{\mathcal L}_n^{(2)})
\otimes \varepsilon_n.
\]
\label{conj2}
\end{conjecture}

More generally, it would be interesting to know an answer to the following
question.

\begin{question} What is the $S_n$-module structure of $\widetilde{H}_{3k-4}
(\nmkn)$ ?
\end{question}  

The complex $\bnmkrs$ also has homology concentrated in
a single dimension. This time the complex is invariant under the natural 
action of $S_r \times S_s$ . Experiments in small cases and
the rank of the homology groups given by Theorem \ref{main-bipartite} suggest
the following conjecture.

\begin{conjecture} As an $(S_r \times S_s)$-module the reduced 
homology group
$\widetilde{H}_{2k-3}(\bnmkrs)$ is isomorphic to the tensor product
$H(r,k) \otimes H(s,k)$, where $H(n,l)$ is the irreducible $S_n$-module
corresponding to the hook Young diagram with $l$ rows and $n-l+1$ columns. 
\end{conjecture}

One can consider the complex $\nmkn$ as the complex of all $n$-partite
graphs on $n$ vertices that do not have a matching of size $k$.
Let $\tnpm$ be the complex of $t$-partite graphs on a $t$-partition with
parts of size ${\mathsf s_1}, \ldots, {\mathsf s_t}$ with no matching of
size $k$. Computational evidence indicates that in general this complex has 
homology concentrated in a single dimension.

\begin{question}
What is the homotopy type of $\tnpm$ ?
Is it homotopic to a wedge of spheres, all of the same dimension ?
\end{question}

The remainder of this paper is organized as follows.  In Section
\ref{gt} we introduce our graph theoretic notation and discuss the
Gallai-Edmonds structure theorem.  In Section \ref{tot} we
introduce certain graphs, called forests of triangles, which will
play a prominent role in our arguments. Namely, these graphs will
represent critical cells of discrete Morse functions on some of
our complexes. Discrete Morse theory is discussed briefly in
Section \ref{dmt}. In Section \ref{prf-a} we give the proofs of
Theorems \ref{fcth} and \ref{bfcth}.  Then in Section \ref{prf-b} we
deduce Theorems \ref{main-general} and \ref{main-bipartite} from
these theorems.

\section{Preliminaries} \label{pre}

\subsection{Graph theory} \label{gt}

By a graph we always mean $G=(V(G),E(G))$ where $V=V(G)$ is a
finite vertex set and $E=E(G) \subseteq \binom{V}{2}$ is the edge
set of $G$.  An edge $\{x,y\} \in E$ will often be denoted by
$xy$.  For $X \subseteq V$, the subgraph of $G$ induced on $X$
will be denoted by $G|_X$, so $G|_X=(X,E \cap \binom{X}{2})$.  For
$v \in V$, $G-v$ will denote $G|_{V -\{v\}}$.  For
distinct $v,w \in V$, $G+vw$ will denote the graph $(V,E \cup
\{vw\})$ obtained by adding $vw$ to $E(G)$ and $G-vw$ will denote
the graph $(V,E \setminus \{vw\})$ obtained by removing $vw$ from
$E(G)$ (so $G+vw=G$ if $vw \in E$ and $G-vw=G$ if $vw \not\in E$).
A matching in $G$ is a subset $M$ of $E(G)$ such that each
$v \in V$ is contained in at most one $e \in M$.  We say a
matching $M$ covers $v \in V$ if some $e \in M$ contains
$v$. A matching $M$ is perfect if every $v \in V$ is covered
by $M$.  Thus we can reformulate the definition of being factor critical 
from Section \ref{intro} as follows. If $|V|$ is odd, we call $G$ factor
critical if for
each $v \in V$, $G-v$ contains a perfect matching.  A maximum
matching in $G$ is one which contains at least as many edges as
every other matching in $G$.  The size of each maximum matching of
$G$ will be denoted by $\nu(G)$.  For $v \in V(G)$, the 
neighborhood $N_G(v)$ is the set of all $w \in V(G)$ such that
$vw \in E(G)$.

For a graph $G=(V,E)$, set
\begin{itemize}
\item $D=D(G):=\{v \in V:$ Some maximum matching of $G$ does
not cover $v\}$, 
\item $A=A(G):=\{v \in V - D~|~N_G(v)
\cap D \neq \emptyset\}$, and 
\item $C=C(G):=V \setminus (A \cup
D)$.
\end{itemize}

\begin{theorem}[Gallai-Edmonds structure theorem, see
{\cite{LoPl}}] For every graph $G$, the following conditions
hold.
\begin{enumerate}
\item Each connected component of $G|_D$ is factor critical. 
\item
Each maximum matching of $G$ consists of
\begin{itemize}
\item a perfect matching on $G|_C$, 
\item for each $a \in A$ an
edge $ad_a$, such that $d_a \in D$ and for distinct $a,b \in A$,
$d_a$ and $d_b$ are in distinct connected components of $G|_D$, and
\item a matching of size $\frac{|V(X)|-1}{2}$ on each connected
component $X$ of $G|_D$.
\end{itemize}
\item Let $\con(G)$ be the number of connected components of $G|_D$.
The number of vertices of $G$ not covered by a given maximum
matching is $\con(G)-|A|$.  In other words,
\[
\nu(G)=\frac{1}{2}\left( |V|-\con(G)+|A| \right).
\]
\end{enumerate}
\label{ge}
\end{theorem}

\subsection{Trees of triangles} \label{tot}

\begin{definition}
Let $V=\{v_1,v_2,\ldots,v_n\} \subseteq \nn$ with $v_i<v_{i+1}$
for all $i<n$.  A connected graph $G$ on vertex set $V$ is a 
tree of triangles if either $|V|=1$ or
\begin{itemize}
\item $v_1v_2 \in E(G)$, 
\item there exists a unique $m \in
\{3,\ldots,n\}$ such that $v_1v_m \in E(G)$ and $v_2v_m \in E(G)$,
and 
\item the graph obtained from $G$ by removing edges
$v_1v_2,v_1v_m,v_2v_m$ has three connected components, each of
which is a tree of triangles.
\end{itemize}
A forest of triangles on $V$ is a graph $F$ on $V$ such that each
connected component of $F$ is a tree of triangles on its vertex
set.
\end{definition}

Induction on $n$ shows that if there is a tree of triangles with $n$ vertices
then $n$ is odd.  For odd $n \in \nn$ we write $n!!$ for the product of all
odd $j \in \nn$ such that $j \leq n$.  By convention, $(-1)!!=1$.

\begin{proposition}
Let $V \subseteq \nn$ with $|V|=2k-1$ for some $k \in \nn$.  Then
the number $\ttk$ of trees of triangles on vertex set $V$ is
$(2k-3)!!^2$. Each tree of triangles on $V$ is factor critical and
has $3k-3$ edges. \label{ntt}
\end{proposition}

\begin{proof}
We proceed by induction on $k$, the case $k=1$ being trivial.  Assume $k>1$.
We may assume that $V=[2k-1]$.  Note that
\begin{equation} \label{ezid}
(2k-3)!!=\frac{(2k-2)!}{(k-1)!2^{k-1}}.
\end{equation}
In each tree of triangles $G$ on vertex set $V$, there is a unique
$v \in \{3,\ldots,2k-1\}$ such that $12,1v,2v \in E(G)$.  We may
choose $v$ in $2k-3$ ways, and having chosen $v$, the removal of
edges $12,1v,2v$ leaves a graph with three connected components
$D_1,D_2,D_v$, such that $i \in D_i$ and the subgraph of $G$
induced on each $D_i$ is a tree of triangles.  It now follows that
$|E(G)|=3k-3$ as claimed, and that if we let $\Pi^0(v)$ be the set
of all ordered partitions of $\{3,\ldots,2k-1\} \setminus \{v\}$
into three parts of even size, then we have

$$
\begin{array}{rcl} 
\ttk & = & (2k-3)\displaystyle{\sum_{(I,J,L) \in \Pi^0(v)}}
\tt(|I|+1)\tt(|J|+1)\tt(|L|+1) \\
     &    & \\
     & =  & (2k-3) \\
     & & \displaystyle{\sum_{i+j+l=k-2}}
\binom{2k-4}{2i,2j,2l} \tt(2i+1)\tt(2j+1)\tt(2l+1) \\ 
     &    & \\
     & =&   \frac{(2k-3)!}{2^{2k-4}} \displaystyle{\sum_{i+j+l=k-2}}
\binom{2i}{i}\binom{2j}{j}\binom{2l}{l},
\end{array}
$$
the last equality following from the inductive hypothesis and equation
(\ref{ezid}). Our formula for the number of trees of triangles now
follows from the facts (easily verified using Taylor's theorem)
that
\[
(1-4x)^{-1/2}=\sum_{i \geq 0}\binom{2i}{i}x^i
\]
and
\[
(1-4x)^{-3/2}=\sum_{j \geq 0}\binom{2j}{j}(2j+1)x^j.
\]
It remains to show that every tree of triangles $G$ is factor
critical.  For any $x \in V$, take $i \in \{1,2,v\}$ with $x \in
D_i$ (where the $D_i$ are as defined above).  Then $D_i-x$ has a
perfect matching by inductive hypothesis, as do both $D_j-j$ for
$j \neq i$.  A perfect matching on $G-x$ is obtained by taking
these three perfect matchings and the edge $jk$, where
$\{i,j,k\}=\{1,2,v\}$.
\end{proof}

\subsection{Discrete Morse theory} \label{dmt}

Here we give a brief summary of some results from
discrete Morse theory for regular cell complexes as developed by
Forman (see \cite{Fo}).
We will only consider cell complexes which are obtained by
starting with a simplicial complex $\Delta$ and taking its
quotient by a (possibly empty) subcomplex $\Gamma$. 
A simplicial
complex $\Delta$ is realized as a CW-complex with one $k$-cell for
each of its $k$-dimensional faces and gluing maps determined by
its face structure. 
Note that we do not distinguish between an (abstract) simplicial
complex and its geometric realization. 

Let $\Gamma$ be any simplicial subcomplex of
$\Delta$. If $\Gamma$ is empty then we define
$\Delta/\Gamma=\Delta$. Otherwise, $\Delta/\Gamma$ is a cell
complex with one point $p_0$ and one $k$-cell $\sigma$ for each
face $\sigma$ of $\Delta$ which is not contained in $\Gamma$. If
$\sigma,\tau$ are two such faces, then any portion of the boundary
$\partial\sigma$ of $\sigma$ which was glued to $\tau$ in $\Delta$ remains
glued to
$\tau$ in the same manner, while any portion of $\partial\sigma$ which was
glued in $\Delta$ to a face of $\Gamma$ is now glued to $p_0$.

In order to formulate discrete Morse theory we need to introduce
some concepts from the theory of directed graphs. 
A directed graph $\D=(V,\A)$ on vertex set $V$ is determined by its
set of arcs (i.e. directed edges) $\A \subseteq V \times V$.
If $W \subseteq V$ then we denote by $\D|_W = (W,\A \cap W \times W)$
the digraph induced by $\D$ on $W$.
Let $\D=(V,\A)$ be any directed graph (digraph) which has no
directed cycle, in particular no $\D$ has no loops.  
A Morse matching on $\D$ is a subset $\M$ of
$\A$ such that
\begin{itemize}
\item[(M1)] each vertex of $\D$ is the head or tail of at most one arc in
$\M$,
and 
\item[(M2)] the digraph $\D_{\M}$ obtained from $\D$ by reversing the
direction of each arc in $\M$ has no directed cycle.
\end{itemize}
A critical cell of a Mose matching $\M$ on $\D$ is
a vertex of $\D$ which is not the head or tail of an arc in $\M$.

If $\P$ is a partially ordered set (poset),
$\D(\P)=(V(\P),\A(\P))$ is the digraph obtained by directed each
edge in the Hasse diagram of $\P$ downwards, that is, the arcs in
$\A(\P)$ are pairs $(x,y)$ where $x$ covers $y$ in $\P$.  If
$\Delta$ is a simplicial complex, let $\D=\D(\Delta)$ be the
directed graph with one vertex for each face of $\Delta$
(including the empty face) and an arc $(\sigma,\tau)$ (directed
from $\sigma$ to $\tau$) whenever $\tau$ is a codimension one face
of $\sigma$.  That is, $\D(\Delta)=\D(\P\Delta)$, where $\P\Delta$
is the poset of faces of $\Delta$.

\begin{theorem}[Forman]
Let $\Delta$ be a simplicial complex and let $\Gamma$ be a
subcomplex of $\Delta$.  Let $\M$ be a Morse matching on
$\D(\Delta)$ such that every face of $\Gamma$ is a critical cell
of $\M$. If $\Gamma$ is the empty complex, assume that the empty
face is not a critical cell of $\M$.  Then the quotient complex
$\Delta/\Gamma$ has the homotopy type of a CW-complex with one
vertex $p$ along with one $k$-cell for each $k$-dimensional
critical cell of $\M$.  In particular, if the critical cells of
$\M$ all have the same dimension $k$ then the given complex has
the homotopy type of a wedge of $k$-dimensional spheres, one sphere for each
critical cell. \label{for}
\end{theorem}

The previous theorem is not stated explicitly in the work of Forman,
but it is an immediate consequence of Forman's theory and
the fact that the
incidence of two adjacent cells in $\Delta/\Gamma$ none of which 
is the distinguished cell $p_0$ is regular. 

The next two elementary results are useful in confirming that a
set $\M$ of arcs is a Morse matching.

\begin{lemma}[{\bf Cluster Lemma} - {\cite[Lemma 2]{Jo}}]
Let $\P_1,\ldots,\P_r$ be pairwise disjoint, order convex
subposets of $\P$. For each $i \in [r]$, let $\M^i$ be an acyclic
matching on $D(\P_i)$. Define a relation on the $\P_{i}$ by
$\P_{i}\le_{c} \P_{j}$ if there exist $x\in \P_i$ and $y\in \P_j$
such that $x \leq y$. Assume that the $\P_i$ satisfy the condition
\begin{itemize}
\item[(${\mathcal P}$)]  The relation $\le_{c}$ defines a partial
order on the $\P_{i}$'s.
\end{itemize}
Then
\[
\M:=\bigcup_{i=1}^{r}\M^i
\]
is an acyclic matching on $\D(\P)$.
\label{cluster}
\end{lemma}

\begin{lemma}[{\bf Cycle Lemma} - {\cite[Proposition 3.1]{Sh}}]
Let $\P$ be an order convex subposet of the face poset of a
simplicial complex $\Sigma$ and assume that $\M \subseteq \A(\P)$
satisfies condition (M1).  Then every directed cycle in $\D_{\M}(\P)$
is of the form
$\sigma_1,\tau_1,\sigma_2,\tau_2,\ldots,\sigma_{r-1},\tau_{r-1},
\sigma_r=\sigma_1$, where
\begin{enumerate}
\item $r \geq 3$, 
\item for each $i \in [r-1]$, there is some $x_i
\in \tau_i$ such that $\tau_i=\sigma_i \cup \{ x_i \}$ and
$(\tau_i,\sigma_i) \in \M$, 
\item for each $i \in [r-1]$, there is
some $y_i \in \tau_i$ such that $\sigma_{i+1}=\tau_i \setminus \{
y_i \}$, and 
\item the multisets $\{ x_i~|~i \in [r] \}$ and $\{
y_i~|~i \in [r] \}$ are equal.
\end{enumerate}
\label{js1}
\end{lemma}

\section{Proofs of Theorems \ref{fcth}, \ref{bfcth}}
\label{prf-a}

\subsection{The complex of factor critical graphs} \label{fcg}

Here we prove Theorem \ref{fcth}, which follows immediately from
Proposition \ref{ntt},  Theorem \ref{for} and the following
result.  Recall that $\Sigma(n)$ is the simplex on vertex set
$\binom{[n]}{2}$.

\begin{lemma} 
Let $n \in \nn$ be odd.  Then there exists a Morse matching $\M$
in $\D(\Sigma(n))$ whose critical cells are the trees of triangles
and the not-factor-critical graphs on vertex set $[n]$. \label{mf}
\end{lemma}

As mentioned in the introduction, the use of the Gallai-Edmonds
structure theorem as a fundamental ingredient in our proof of
Lemma \ref{mf} distinguishes this proof from those of previous
results which use discrete Morse theory in examining monotone
graph properties. The rest of the proof is similar in form and
spirit to many of the previous proofs, see
\cite{BBLSW,Jo,LiSh,Sh}.  We will not give as many details as
were given in these proofs.  In particular, at several junctures
we will leave  it to the reader to confirm that we have used
Cluster Lemma \ref{cluster} appropriately or that one can use Cycle 
Lemma \ref{js1} to show that a given matching is actually a Morse
matching.

\begin{proof}
We prove Lemma \ref{mf} by induction on $n$, the case $n=1$ being
trivial.  So, assume $n>1$.  We will construct our Morse matching
$\M$ in several steps.
\\
\\
\noindent {\bf Step 0:} We begin with $\M$ empty and then add to $\M$ the
set $\M^0$ all arcs $(G+12,G)$ where $G$ is factor critical but
$12 \not\in E(G)$. Using Lemma \ref{js1}, it is easy to confirm
that these arcs form a matching and that their reversal leaves an
acyclic digraph. Certainly every graph which is covered by an arc
in $\M^0$ is factor critical.  The factor critical graphs which
are not covered by any arc in $\M^0$ are those factor critical
graphs $G$ such that $12 \in E(G)$ and $G-12$ is not factor
critical. Let $\C^0$ be the set of all such graphs, 

In the following we collect properties of graphs in $\C^0$.
Let $G \in \C^0$ (so $G$ is factor critical but $G-12$ is not). 
Set
\[
\gp:=G-12.
\]
We consider the Gallai-Edmonds decomposition of $\gp$.  Since $G$
is factor critical, we see that $G-1 = \gp-1$
and $G-2 = \gp-2$ both have perfect matchings.  This gives

\begin{equation} \label{sizematching}
\nu(\gp)=\frac{n-1}{2}
\end{equation}

and

\begin{equation} \label{12indgp}
1,2 \in D(\gp).
\end{equation}

\noindent {\sf Claim 1:} The vertices $1$ and $2$ lie in different
connected components of
$\gp|_{D(\gp)}$.

\noindent 
$\triangleleft$ 
{\sf Proof of Claim 1:}
Assume for contradiction that $1,2$ lie in the same component
$X=(V(X),E(X))$ of $\gp|_{D(\gp)}$.  
Since $\gp$ is not factor critical it follows that there is a vertex $v$
such that $\gp - v$ does not have a perfect matching. Since $\nu(\gp) =
\frac{n-1}{2}$, we see that $v$ is not in  $D(\gp)$. Consequently,
$v \in  A(\gp) \cup
C(\gp)$.  For any such $v$, there is a perfect matching $K$ in $G-v$, and
since $v \not\in
D(\gp)$, this matching $K$ contains $12$.  By Theorem \ref{ge}(1),
there are oddly many elements of $V(X) \setminus \{1,2\}$ and since all of 
them are contained in an edge from $K$ it
follows that there exist $x \in V(X)$ and $a \in A(\gp)$ such that
$xa \in K$. By Theorem \ref{ge}(3) and $\nu(\gp)=\frac{n-1}{2}$ it follows
that $\con(\gp) = |A|+1$. But now there are $\con(\gp)-1$ components of
$\gp|_{D(\gp)}$ other than $X$, each of which have odd size and
$|A(\gp)|-1<\con(\gp)-1$ elements of $A(\gp) \setminus \{a\}$
remaining to pair with elements of $D(\gp) \setminus V(X)$ in $K$.
Thus by Theorem \ref{ge}(2) $K$ cannot be a perfect matching, giving the 
desired contradiction. $\triangleright$

\noindent {\sf Claim 2:} $A(\gp) \neq \emptyset$.

\noindent 
$\triangleleft$ 
{\sf Proof of Claim 2:} 
By (\ref{sizematching}) and Theorem \ref{ge}(3) it
follows that $\con(\gp) = 1 + |A(\gp)|$. From {\sf Claim 1} we know
that $\con(\gp) \geq 2$ which implies the assertion. $\triangleright$

\medskip

\noindent {\bf Step 1:} 
The set of graphs in $\C^0$ forms an ideal in 
the poset of factor critical graphs, so we can apply Cluster Lemma
\ref{cluster} after describing a Morse matching on the subgraph
$\D^0$ of $\D(\Sigma(n))$ induced on $\C^0$.

For all $G \in \C^0$ such that $|A(\gp)|>1$, let $a(G),b(G)$ be the
two smallest elements of $A(\gp)$.  Define the set $\M^1$ of arcs
in $\D^0$ by
\[
\M^1:=\{(G+a(G)b(G),G)~|~G \in C^0,a(G)b(G) \not\in E(G)\}.
\]
We will see that $\M^1$ is a Morse matching on $\D^0$.  First note
that for any graph $H$ and any distinct $a,b \in A(H)$, the graphs
$H+ab$ and $H-ab$ have the same Gallai-Edmonds decomposition,
since (by Theorem \ref{ge}(2)) no maximum matching in $H$ uses an
edge with endpoints in $A(H)$.  It follows that no vertex of
$\D^0$ is a head or tail of more than one arc in $\M^1$, and it
remains to show that $\D^0_{\M^1}$ has no directed cycle. Assume
for contradiction that such a directed cycle exists.  By Cycle
Lemma \ref{js1}, this cycle has vertices
\[
G_1,H_1,\ldots,G_{r-1},H_{r-1},G_r=G_1,
\]
where
\begin{itemize}
\item $r \geq 3$, 
\item $H_i=G_i-xy$ for some $xy \in E(G_i)$, and
\item $G_{i+1}=H_i+a(\hp_i)b(\hp_i)$.
\end{itemize}

As noted above, $\hp_i$ and $\gp_{i+1}$ have the same
Gallai-Edmonds decomposition for all $i$.  Moreover, since all the
$G_i$ and $H_i$ lie in $\C^0$, we have
\[
\nu(\hp_i)=\nu(\gp_i)=\frac{n-1}{2}
\]
for all $i$. Thus any maximum matching of $\hp_i$ is a maximum matching of
$\gp_i$.
It follows that $D(\gp_i) \subseteq D(\hp_i)$ for
all $i$.  Since $G_r=G_1$, we must have
\[
D(\hp_i)=D(\gp_i)=D(\gp_1)
\]
for all $i$.

\noindent {\sf Claim 3:} $A(\hp_i)=A(\gp_i)$ for all $i$. \\

\noindent {\sf Claim 4:} From {\sf Claim 3} the desired contradiction
follows. \\

\noindent 
$\triangleleft$ 
{\sf Proof of Claim 4:} 
Indeed, from the validity of {\sf Claim 3} 
it follows that $a:=a(\hp_1)=a(\gp_1)$ and $b:=b(\hp_1)=b(\gp_1)$.  
If $\hp_1=\gp_1-ab$ then $\gp_2=\gp_1$, a contradiction.  If $ab \in
E(\hp_1)$ then there is no arc $(K,H_1)$ in $\M^1$, again a
contradiction.  Otherwise, we have $\gp_2=\hp_1+ab$.  But then either
$\hp_2 = \gp_2 - ab = \hp_1$ (a contradiction) or $ab \in \hp_2$, in which case
$A(\hp_2)=A(\gp_2)$ by {\sf Claim 3}. Thus 
there is no arc $(K,H_2)$ in $\M^1$, a contradiction. $\triangleright$ 

\noindent 
$\triangleleft$ 
{\sf Proof of Claim 3:} 
Since $\nu(\hp_1) = \nu(\gp_1)$ it follows
from Theorem \ref{ge}(3) that $\con(\hp_1) - |A(\hp_1)|=\con(\gp_1) -
|A(\gp_1)|$.
Since $D(\hp_1) = D(\gp_1)$ and $\hp_1 = \gp_1 -xy$, we either have 
$\con(\hp_1) = \con(\gp_1)$ or
$\con(\hp_1) = \con(\gp_1)+1$. In the first case $|A(\gp_1)|= |A(\hp_1)|$.
In the second case $xy$ must be an edge whose removal disconnects a connected
component of $\gp|_{D(\gp)}$ and $|A(\gp_1)| = |A(\hp_1)|+1$. But again by
$D(\gp_1) = D(\hp_1)$ this implies that $xy$ is an edge which connects an
element of $D(\gp_1)$ with an element of $A(\gp_1)$, a contradiction.
Thus we may assume $|A(\gp_1)|= |A(\hp_1)|$ and
it is sufficient to show that
$A(\gp_1) \subseteq A(\hp_1)$.  To prove this last fact, it
suffices to show that for each $a \in A(\gp_1)$ we have
$|N_{\gp_1}(a) \cap D(\gp_1)|>1$. Assume for contradiction that
some $a \in A(\gp_1)$ has only one neighbor $d$ in $D(\gp_1)$.
Then $\gp_1-d$ has no perfect matching, as the $\con(\gp_1)-1$
connected components of the subgraph of $\gp_1$ induced on
$D(\gp_1)$ have together only $|A(\gp_1)|-1$ neighbors in
$A(\gp_1)$. $\triangleright$ 

\smallskip 

This completes the proof that $\M^1$ is a Morse
matching on $\D^0$.  As noted above, Cluster Lemma \ref{cluster}
guarantees that $\M^0 \cup \M^1$ is a Morse matching.
\\
\\
\noindent {\bf Step 2:} Let $\C^1 \subseteq \D^0$ be the set of critical
cells of $\M^0 \cup \M^1$. 

\smallskip

\noindent {\sf Claim 5:} $G \in \C^1$ if and only if
\begin{itemize}
\item[(a)] $G$ is factor critical,
\item[(b)] $|A(\gp)|=1$, and
\item[(c)] $\gp|_{D(\gp)}$ has exactly two connected components
$D_1(\gp),D_2(\gp)$
for which $i \in D_i(\gp)$, $i=1,2$ holds.
\end{itemize}

\noindent 
$\triangleleft$ 
{\sf Proof of Claim 5:}
Since all critical cells of $\C^0$ are factor critical the same holds for
$\C^1 \subseteq \C^0$. From {\sf Claim 2} we know that  $|A(\gp)| \geq 1$ 
for all $G \in \C^1 \subseteq \C^0$. Now exactly those $G \in \C^0$ with
$|A(\gp)| \geq 2$ are covered by an edge from $\M^1$. Thus (b) holds for
all $G \in \C^1$. By {\sf Claim 1} the vertices $1$ and $2$ lie in different
connected components $D_1(\gp)$ and $D_2(\gp)$ of $\gp|_{D(\gp)}$. 
By (\ref{sizematching}) and Theorem \ref{ge}(3) it follows that 
$\con(G) = 1+|A(\gp)|$. From (b) we know that
$|A(\gp)| = 1$ and hence $\con(\gp) = 2$. Thus $D_1(\gp)$ and $D_2(\gp)$
are the only connected components of $\gp|_{D(\gp)}$.

We have shown that every element of $G \in \C^1$ satisfies (a)-(c). 
Conversely it is easily checked that no graph satisfying (a)-(c) is
covered by an edge in $\M^0 \cup \M^1$. $\triangleright$

\smallskip

For each $a \in [n] - \{1,2\}$ and each ordered partition $(X,Y,Z)$ 
of $[n] \setminus \{1,2,a\}$ into three possibly empty subsets we define
the subset $\C^1[a,(X,Y,Z)] \subseteq \C^1$ by
$G \in \C^1[a,(X,Y,Z)]$ if and only if
\begin{itemize}
\item $G \in \C^1$, 
\item $A(\gp)=\{a\}$,
\item $V(D_1(\gp))=X \cup \{1\}$,
\item $V(D_2(\gp))=Y \cup \{2\}$, and
\item $C(\gp)=Z$.
\end{itemize}

It follows from {\sf Claim 5} that the $\C^1[(a,X,Y,Z)]$ actually
partition $\C^1$.  The Cluster Lemma \ref{cluster} applies to this partition, 
and we will define a Morse matching $\M[a,(X,Y,Z)]$ on each 
$\C^1[(a,X,Y,Z)]$.

Fix $a,(X,Y,Z)$ and let $\C=\C^1[a,(X,Y,Z)]$.  We will show that
there is a Morse matching $\M[a,(X,Y,Z)]$ on the subgraph $\D|_\C$
of $\D^0$ induced on $\C$ whose critical cells are exactly those
trees of triangles $G$ such that
\begin{itemize}
\item $12,1a,2a \in E(G)$ and
\item the connected components of $G-\{12,1a,2a\}$ are $X \cup \{1\}$, $Y
\cup \{2\}$
and $Z \cup \{a\}$.
\end{itemize}
Once this is done, the proof of our lemma is completed by applying Cluster
Lemma
\ref{cluster} to $\M^0$, $\M^1$ and all of the $\M[a,(X,Y,Z)]$.  Set
\[
\M(1):=\{(G+1a,G):G \in \C,1a \not\in E(G)\}.
\]
It is straightforward to show that $\M(1)$ is a Morse matching on $\D|_\C$
whose critical
cells are those $G \in \C$ such that $N_G(a) \cap X=\emptyset$.  Let
$\C(1)$ be the
set of all such $G$, and define
\[
\M(2):=\{(G+2a,G):G \in \C(1),2a \not\in E(G)\}.
\]
It is again straightforward to show (using Cluster Lemma \ref{cluster})
that $\M(1) \cup
\M(2)$ is a Morse matching on $\D|_\C$ whose critical cells are those $G
\in \C(1)$ such that
$N_G(a) \cap Y=\emptyset$.  Let $\C(2)$ be the set of all such critical cells.

We claim now that if $G \in \C(2)$ then $G|_{Z \cup \{a\}}$ is factor
critical.  In other
words, $\C(2)$ consists of all graphs $G \in \C$ such that
\begin{itemize}
\item $12,1a,2a \in E(G)$,
\item $G-\{12,1a,2a\}$ has three connected components $D_1=X \cup \{1\}$,
$D_2=Y \cup
\{2\}$, $D_a=Z \cup \{a\}$, and
\item the subgraph of $G$ induced on each of these three components is
factor critical.
\end{itemize}
If this claim holds then we can use the inductive hypothesis (and
Cluster Lemma \ref{cluster} one more time) to produce the desired
Morse matching on $\C$ and our lemma follows.

So, let $G \in \C(2)$. Note that since $G|_Z=G|_{C(\gp)}$ contains a
perfect matching by
definition, it remains to show that if $z \in Z$ then the subgraph of $G$
induced
on $(Z \setminus \{z\}) \cup \{a\}$ contains a perfect matching.  We know
that $G-z$
contains a perfect matching but that $\gp-z$ contains no perfect matching.
Therefore,
any perfect matching $K$ in $G-z$ includes the edge $12$.  Since both $|X|$
and $|Y|$ are
even and the connected components of $\gp|_{D(\gp)}$ are $D_1(\gp)$ and
$D_1(\gp)$,
$K$ cannot contain any edge $av$ with $v \in D(\gp)$.  Thus $K$ consists of
$12$, a
perfect matching on $X$, a perfect matching on $Y$ and the desired perfect
matching on
$(Z \setminus \{z\}) \cup \{a\}$ and we are done.
\end{proof}

\subsection{The complex of $q$-factor critical bipartite graphs} \label{bfcg}

The following lemma immediately implies Theorem \ref{bfcth}.

\begin{lemma} For $0 \leq q <s$, there is a Morse matching $\M$ on
the digraph $D(\Sigma(q,s))$ whose critical cells are the elements of
$\nbfc(q,s)$
along with $\binom{s-1}{q}$ cells of dimension $2q-1$ from $\bfc(q,s)$.
\end{lemma}
\begin{proof}
We proceed by induction on $q$, the case $q=0$ being trivial. Assume
$q > 0$. As in the proof of Lemma \ref{mf} we construct our Morse
function in several steps.

\noindent {\bf Step 0:} We begin with $\M$ being empty and then
add to $\M$ the set $\M^0$ of all arcs $(G+1\overline{1},G)$ where
$G$ is $q$-factor critical and $1\overline{1} \not\in E(G)$. Then $\M^0$
is a matching, the digraph $D_{\M^0}$ is acyclic and the set $\C^0$ of
graphs not covered by $\M^0$ consists of $\nbfc(q,s)$ along with the elements
$G \in \bfc(q,s)$ such that $1\overline{1} \in E(G)$ and $G -
1\overline{1} \not\in \bfc(q,s)$.

For $G \in \C^0 \cap \bfc(q,s)$, set $\gp := G - 1\overline{1}$ and consider
the Gallai-Edmonds decompositions of the vertex set of $\gp$ into $A', C'$ and
$D'$. We prove the following claims:

\noindent {\sf Claim 1:} $\overline{1} \in D'$ and $\nu(\gp) = q$.

\noindent $\triangleleft$ {\sf Proof:} Since $G \in \bfc(q,s)$, we have 
$\nu(G - \overline{1}) = q$. The assertion follows from $G - \overline{1} =
\gp - \overline{1}$. $\triangleright$

\noindent {\sf Claim 2:} $D' \subseteq [\overline{s}]$ and therefore $A'
\subseteq [q]$.

\noindent $\triangleleft$ {\sf Proof:} Since by {\sf Claim 1} we have
$\nu(\gp) = q$  
and $\gp$ is bipartite, we must have $D' \subseteq [\overline{s}]$. The second
part of the claim then follows again from the fact that $\gp$ is bipartite.
$\triangleright$

\noindent {\sf Claim 3:} $1 \in C'$.

\noindent $\triangleleft$ {\sf Proof:} 
By {\sf Claim 2} it suffices to show that $1 \not\in A'$. 
Assume for contradiction that $1 \in A'$.  
For each $y \in D'$ we have $N_{G'}(y) \subseteq A'$ and the assumption
$1\in A'$ implies $N_G(y) \subseteq A'$. 
Since $\nu(G')=q$ by {\sf Claim 1} and
$G' \not\in \bfc(q,s)$, we see that $[\overline{s}] \neq D'$.  Thus there
is some $\overline{x} \in [\overline{s}] \cap C'$.  Let $p=|[\overline{s}]
\cap C'|$.  Note that $p=|[q] \cap C'|$, since $C'$ contains a perfect
matching.
{}From $\nu(G')=q$ we infer $|A'|=q-p$ and $|D'|=s-p$.  Let $K$ be a
matching of
size $q$ in $G-\overline{x}$. ($K$ exists since $G \in \bfc(q,s)$.)  
At most $s-q$ elements of $D'$ are not covered by any edge of $K$.
Since $D'$ has $s-p$ elements, all in $[\bar{s}]$, there must be at least
$q-p$ edges of $K$ with one endpoint in $D'$ and one endpoint in $A'$.  
Since $|A'|=q-p$, every element of $A'$ is the
endpoint of one of the edges just mentioned.  This means that there is no
edge in $K$ with one endpoint in $C'$ and one in $A'$, and certainly there is
no edge in $K$ with one endpoint in $C'$ and one in $D'$.  Thus the
remaining $p$ edges in $K$ have both endpoints in $C'\setminus
\{\overline{x}\}$.  But this is impossible since $|(C' \setminus
\{\overline{x}\}) \cap [\overline{s}]|<p$.
$\triangleright$

\noindent {\bf Step 1:} Set $\B^0 :=  \{ G \in \C^0 \cap \bfc(q,s)~|~A' 
\neq \emptyset\}$, carrying forward the definitions of
$A', C'$ and $D'$ from {\bf Step 0}. For $G \in \B^0$, let $a$
(resp. $\overline{c}$) be the smallest elements of $A'$ (resp. 
$C' \cap [\overline{s}]$). Note that by {\sf Claim 3} we have
$a \neq 1$. Define $\M^1 := \{ (G,G-a\overline{c})~|~
G \in \B^0, a\overline{c} \in E(G)\}$. 
 
\noindent {\sf Claim 4:} For $G \in \B^0$ such that 
$a\overline{c} \in E(G)$ the graphs $\gp$ and $\gp -
a\overline{c}$ have the same Gallai-Edmonds decomposition.

\noindent $\triangleleft$ {\sf Proof:} By Theorem \ref{ge}(2) the 
edge $a \overline{c}$ is not included in 
any maximum matching of $\gp$. $\triangleright$
   
It follows immediately from {\sf Claim 4} that $\M^1$ is a matching.

\noindent {\sf Claim 5:} For each $G \in \B^0$, we have $G- a\overline{c}
\in \B^0$.  

\noindent $\triangleleft$ {\sf Proof:} 
By {\sf Claim 4}, it suffices to show that $G - a\overline{c} \in
\C^0 \cap \bfc(q,s)$. Since $G \in \C^0$, it suffices to show that
$G-a\overline{c}
\in \bfc(q,s)$.  Let $\overline{x} \in [\overline{s}]$.  If $\overline{x}
\in D'$
then there is a matching of size $q$ in $G'-\overline{x}$ which does not
contain
$a\overline{c}$ by Theorem \ref{ge}.  Say $\overline{x} \in C'$.  There is
a matching
$K$ of size $q$ in $G-\overline{x}$, and every such matching contains the edge
$1\overline{1}$.  As above, let $p=|[\overline{s}] \cap C'|=|[q] \cap C'|$,
so $|A'|=q-p$
and $|D'|=s-p$.  Let $k$ be the number of edges of $K$ which contain an
element $y$ of $A'$
and an element $\overline{z}$ of $C'$.  Note that since $y \in [q]$, we
must have 
$\overline{z} \in [\overline{s}]$.  

The only edge in $G$ between $C'\cap[q]$ and $D'$ is $1\overline{1}$. Thus
$p-1-(p-k-1)=k$ elements of $C'\cap[q]$ are not covered by $K$. Hence $k=0$
and in particular, $a\overline{c} \not\in K$ as desired.
$\triangleright$

\noindent {\bf Step 2:} Let $\C^1$ be the set of critical points of $\M^0 \cup
\M^1$. Then $\C^1$ consists of $\nbfc(q,s)$ and those $G \in \bfc(q,s)$
such that
$\gp \not\in \bfc(q,s)$ and $A'= \emptyset$. By Theorem \ref{ge}(3) it then
follows that $|D'| = s-q$.   

Note that $\overline{1} \in D'$. For each $X \subseteq [\overline{s}]
- \{ \overline{1} \}$ such that $|X| = s-q-1$, set $\C^1[X] := \{ G \in \C^1 
\cap \bfc(q,s)~ |~ D' = X \cup \{ \overline{1}\}\}$. Note that there are
$\binom{s-1}{s-q-1} = \binom{s-1}{q}$ choices for $X$. We can then apply the
Cluster Lemma \ref{cluster} to the decomposition of $\C^1 $ into
$\nbfc(q,s)$ and the sets $\C^1[X]$. 
 
Fix one such $X$. Let $G \in \C^1[X]$. Then the $s-q$ elements of $D'$ 
are by Remark \ref{critical-bipartite} isolated in $\gp$. 

\noindent {\sf Claim 6:} $G|_{C'} - 1$ is $(q-1)$-factor critical. 

\noindent $\triangleleft$ {\sf Proof:} The claim follows from the fact 
that, if $\overline{x} \in C' \cap [\overline{s}]$ then any
matching of size $q$ in $G - \overline{x}$ consists of the edge 
$1\overline{1}$ along with $q-1$ edges in $C' - \overline{x}$. 
$\triangleright$ 

Write $N_{C'}(1)$ for $N_G(1) \cap C'$

\noindent {\sf Claim 7:} $N_{C'}(1) \neq \emptyset$.

\noindent $\triangleleft$ {\sf Proof:} Otherwise $\nu(G-\overline{1})=q-1$.
$\triangleright$

\noindent {\bf Step 3:}
For fixed $X$ and $G \in \C^1[X]$, let $\overline{c}=\overline{c}(G)$ be the
smallest element of $C' \cap [\overline{s}]$.  Define
\[
\M^2[X]:=\{(G+1\overline{c},G):G \in \C^1[X], \overline{c}(G) \not\in
N_G(1) \}.
\]
It is straightforward to show that $\M^2[X]$ is a Morse matching on
$\D(\Sigma(q,s))|_{\C^1[X]}$.  Let $\C^2[X]$ be the set of critical points of 
$\M^2[X]$.

\noindent {\sf Claim 8:} $G \in \C^2[X]$ if and only if
\begin{itemize} \item Every element of $D'=X \cup \{\overline{1}\}$ is
isolated
in $G$, \item $G|_{C'}-1$ is $(q-1)$-factor critical, and \item 
$N_{C'}(1)=\{\overline{c}(G)\}$.
\end{itemize}

\noindent $\triangleleft$ {\sf Proof:} We have already seen that the first
two conditions
are necessary.  To show that the third condition is also necessary, it
suffices to show
that if $G$ satisfies the first two conditions and $\overline{c}(G) \neq
\overline{y} \in
N_{C'}(1)$ then $G-1\overline{c}(G) \in \C^1[X]$.  Equivalently, we must
show that 
$G-1\overline{c}(G) \in \bfc(q,s)$ and that $D(\gp-1\overline{c}(G))=D'$.

Let $\overline{z} \in [\overline{s}]$.  If $\overline{z} \in D'$ then a
matching of size $q$
in $\gp-1\overline{c}(G)-\overline{z}$ is obtained by taking a matching of
size $q-1$ in 
$G|_{C' \setminus \{1,\overline{y}\}}$ along with $1\overline{y}$.  Thus
$\overline{z} \in 
D(\gp-1\overline{c})$.  Say $\overline{z} \in C'$.  A matching of size $q$ in 
$G-1\overline{c}-\overline{z}$ is obtained by taking a matching of size
$q-1$ in 
$G|_{C'\setminus \{1,\overline{z}\}}$ along with $1\overline{1}$.
Moreover, we have
\[
\nu(\gp-1\overline{c}-\overline{z}) \leq \nu(\gp-\overline{z})<q,
\]
so $\overline{z} \not\in D(\gp-1\overline{c})$.

To show that the three conditions are sufficient, it suffices to show that
$G \in \bfc(q,s)$
but $\gp \not\in \bfc(q,s)$.  If $1$ has no neighbor in $C'$ other than
$\overline{c}(G)$ and
$A'=\emptyset$ then $\nu(\gp-\overline{c})<q$, so $\gp \not\in \bfc(q,s)$.
On the other hand,
if $G|_{C'}-1$ is $(q-1)$-factor critical then for each $\overline{z} \in
[\overline{s}]$ one
obtains a matching of size $q$ in $G-\overline{z}$ by taking a matching of
size $q-1$ in
$G|_{C' \setminus \{1,\overline{z}\}}$ along with $1\overline{1}$.
$\triangleright$

We see now that $\D(\C^2[X])$ is isomorphic to $\D(\bfc(q-1,q))$.  thus by
our inductive 
hypothesis there is a matching $\M^3[X]$ on $\D(\C^2[X])$ whose unique
critical cell is a
graph $G$ with $2q$ edges ($2q-2$ of them having one endpoint in $[q]
\setminus \{1\}$ and
one in $[\overline{s}] \setminus (X \cup \{\overline{1}\})$ and the other
two being 
$1 \overline{1}$ and $1\overline{c}(G)$).

Using Cluster Lemma \ref{cluster}, we see that $\M^0 \cup \M^1 \cup
\Big\{M^2[X] \cup M^3[X]:
X \in \binom{[\overline{s}] \setminus \{\overline{1}\}}{s-q-1} \Big\}$ is
the desired Morse matching.

\end{proof}

\section{Proofs of Theorems \ref{main-general}, \ref{main-bipartite}}
\label{prf-b}

\subsection{The complexes $\nmkn$} \label{nmkn} 

We now prove Theorem \ref{main-general}

\begin{proof}
We proceed in two steps:

\noindent {\bf Step 0:} 
We define a Morse matching $\M$ on $\nmkn$ such that the graphs 
$G$ corresponding to the critical cells are exactly those which satisfy:
\begin{itemize}
\item $N_G(n)=\emptyset$
\item For each $1 \leq v<n$, $G+vn$ contains a matching of size $k$.
\end{itemize}

For $1 \leq v \leq n-1$, define $\M(i)$ on $\nmkn$ recursively as follows.
\begin{itemize}
\item $\M(1) :=\{(G,G-1n)~|~G \in \nmkn,1n \in E(G)\}$. 
\item For $1<v<n$, let $\C(v)$ consist of those $G \in \nmkn$ such that
no arc in $\bigcup_{w<v}\M(w)$ has $G$ as its head or tail. Then
$\M(v):=\{(G,G-vn)~|~G \in \C(v),vn \in E(G)\}$.
\end{itemize}

Note that the presence or absence in a graph $G$ of an edge $wn$
($w<v$) has no effect on the existence of a matching of size $k$
in $G$ which includes the edge $vn$.  Thus it follows from Cluster
Lemma \ref{cluster} that $\M := \bigcup_{v=1}^{n-1}\M(v)$ is a
Morse matching on $\nmkn$.  

The critical cells of this Morse
matching are those $G \in \nmkn$ such that $N_G(n)=\emptyset$ and,
for $1 \leq v<n$, $G+vn$ contains a matching of size $k$.  
For each such $G$ and for $1 \leq v<n$, we have $\nu(G-v)=k-1$.  

\noindent {\bf Step 1:} 
Let $G$ be a graph corresponding to a critical cell of the Morse
matching $\M$  from {\bf Step 0}. Let $\gm :=G-n$.  
Since for $1 \leq v<n$, $G+vn$ contains a matching of size $k$
we have $\nu(\gm-v)=k-1$ and $D(\gm) = [n-1] = V(\gm)$.

Thus by Theorem \ref{ge}(3) each connected
component of $\gm = \gm|_{D(\gm)}$ is factor critical.  
Since $\nu(G)=k-1$ and each component $X$ of $\gm$ satisfies
$\nu(\gm |_{X})=\frac{|X|-1}{2}$, we have

\[
c:=\con(\gm)=n-2k+1.
\]

Now for each partition $\tau$ of $[n-1]$ into $c$ parts, let
$\D^\tau$ be the subdigraph of $\D(\nmkn)$ induced on those
critical cells $G$ described above such that the connected
components of $\gm$ determine the partition $\tau$.  Using Lemma
\ref{mf} and Cluster Lemma \ref{cluster}, we can construct a Morse
matching on $\D^\tau$ whose critical cells are those forests of
triangles whose connected components determine $\tau$.  Theorem
\ref{main-general} follows from a final application of Cluster Lemma
\ref{cluster}.
\end{proof}

\subsection{The complexes $\bnmkrs$} \label{bnpkrs}

We now prove Theorem \ref{main-bipartite}.
Analogous to the proof for non-bipartite graphs we have to proceed in two
steps of which {\bf Step 0} is only a slight modification of {\bf Step 0}
from the proof of Theorem \ref{main-general}.

\begin{proof} 

\noindent {\bf Step 0:} 
We inductively define a Morse matching $\M$ on $\bnmkrs$.

For $\bar{1} \leq \bar{v} \leq \bar{s}$, define $\M^{\bar{v}}$ on $\bnmkrs$
recursively as follows.
\begin{itemize}
\item $\M^{\overline{1}} :=\{(G,G-r\bar{1})~|~G \in \bnmkrs,r\bar{1} \in
E(G)\}$. 
\item For $\overline{1} <  \bar{v} \leq \overline{s}$, let $\C^{\bar{v}}$
consist of those 
$G \in \bnmkrs$ which are critical cells of
$\bigcup_{\bar{w}<\bar{v}}\M^{\bar{w}}$, and set
\[
\M^{\bar{v}}:=\{(G,G-r\bar{v}):G \in \C^{\bar{v}},r\bar{v} \in E(G)\}. 
\]
\end{itemize}

Set $\M := \bigcup_{\bar{v} = \overline{1}}^{\overline{s}}
\M^{\overline{v}}$. 
With the same arguments as for non-bipartite graphs we see that
$\M$ actually is a 
Morse matching on $\bnmkrs$.  

The graphs
$G$ corresponding
to the critical cells of $\M$ are exactly those which satisfy:
\begin{itemize}
\item $N_G(r)=\emptyset$.
\item For each $\bar{1} \leq \bar{v}<\bar{s}$, $G+r\bar{v}$ contains a
matching of size $k$.
\end{itemize}

\noindent {\bf Step 1:} 
It follows from K\"onig's Theorem (see \cite{LoPl}), which says that the
size of a 
maximal matching equals the size of a minimal vertex cover,
that exactly $k-1$ elements of $[r]$ are not
isolated.  So, we make one of $\binom{r-1}{k-1}$ choices for the set of
nonisolated
vertices in $[r]$ and then assume that $r = k$. In this case, remaining
graphs are those
for which $[r-1]$ can be completely matched into every 
$(s-1)$-subset of $[\overline{s}]$.  Theorem \ref{main-bipartite} now
follows from (the
proof of) Theorem \ref{bfcth} and Cluster Lemma \ref{cluster}.
\end{proof}

\end{document}